# Exponential Observation Error in Boundary Region

ANAS DHEYAB AL-JOUBORY[1*]

1 School of Mathematics and Statistics, Huazhong University of science and technology

Wuhan, Hubei 430074, PR China

*E-mail: anas_dheyab@hust.edu.cn, tel:+8613164622606, ORCID: 0000-0003-1606-7835.

Abstract:

In this paper the exponential perception blunder idea on account of limit locale has been discussed and broke down. For disseminated parameter frameworks of explanatory sort, we demonstrate that, the blunder of state reproduction can be diminishes by exponentially perception.

## 1. Introduction

In current numerical control hypothesis, perceptibility implies that it is conceivable to recreate remarkably the underlying condition of the dynamic framework from the learning of the information and yield [1]. Thought of provincial discernibleness (reached out by El Jai et al. [2-3]) is of awesome significance in flow look into and roused by numerous applications [4-5,7-8]. The idea of local asymptotic investigation was investigated by Al-Saphory and El Jai [9, 11], and these concepts will be extended to the bundary case by A. Aljoubory and Al-Saphory such that the behavior of the system not in whole domain $\Omega$ but only on particular region $\Gamma$ of the domain $\Omega$ [6, 10]. The purpose of this paper is to talk about the link between the regeonal bundary exponential observation and error (figure 1).

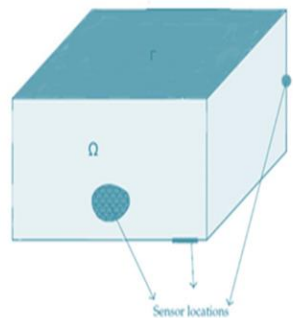

**Figure1. The domain Ω, the boundary region Γ, and the sensors locations**

We suppose that a class of appropriated framework and create distinctive outcomes associated with the different sorts of estimations, and we demonstrate that, the mistake of state remaking can be diminishes by exponentially perception.

This paper is sorted out as takes after: segment one is centered around preliminaries and the issue plan. In the following area, the portrayal thought of provincial limit exponential perception blunder is given.

## 2. Regional boundary exponential observation

In this segment we give initially, the announcement of the issue with the speculation of consdered system, and afterward the idea of provincial limit exponential perceptibility is clarified, and we give a hypothesis which gives the approach watched the present state $z(\mu, t)$ of the original system (1) exponentially.

### 2.1 Preliminaries

The considered framework is portrayed by the explanatory equations;

$$\begin{aligned}
\frac{\partial z}{\partial t}(\mu, t) &= Az(\mu, t) + Bu(t) & \Omega \times (0, \infty) \\
z(\mu, 0) &= z_0(\mu) & \overline{\Omega} \\
\frac{\partial z}{\partial v_A}(\mu, t) &= 0 & \partial\Omega \times (0, \infty)
\end{aligned} \quad (1)$$

Where $\Omega$ is the space when the above framework is characterized as limited open subset of $R^n$ with boundary $\partial\Omega$, $[0,T]$ is the time interval for $T > 0$, $A$ will be a moment arrange direct differential administrator and is self adjoint with reduced resolvant, furthermore, which creates a firmly ceaseless semi-group $(S_A(t))_{t \geq 0}$ on the state space $Z = H^1(\Omega)$ which is one order Sobolv space, $A^*$ will denote the adjoint operator of A. The operators $B \in L(R^p, Z)$ and $C \in L^2(0, T, R^q)$, such that $p$ & $q$ is nmber of controls and estimators. The initial state $z_0(\mu)$ is to be supposed obscure and found inside $H^1(\overline{\Omega})$. The estimations of framework (1) are acquired through inside or limit zone or pointwise sensors which portray the yield work:

$$y(., t) = Cz(\mu, t) \quad (2)$$

Under the above presumption, the framework (1) has a one of a kind arrangement [1-2]:

$$z(\mu, t) = S_A(t)z_0(\mu) + \int_0^t S_A(t-s)Bu(s)ds \tag{3}$$

Presently, we characterize the accompanying administrators:

$$K: Z \to \mathcal{O} \tag{4}$$
$$z \to CS_A(.)z$$

with an adjoint

$$K^*: \mathcal{O} \to Z \tag{5}$$

given by

$$K^*y^* = \int_0^t S_A^*(s)C^*y^*(s)ds \tag{6}$$

We also consider the trace operator of order zero

$$\gamma_0: H^1(\Omega) \to H^{\frac{1}{2}}(\partial\Omega) \tag{7}$$

Which is linear, subjective, and continuous, such that $z_0^\Gamma$ is the restriction of the trace of the initial state $z_0$ to $\Gamma$.

$\gamma_0^*$ denote the adjoint of $\gamma_0$ given by

$$\gamma_0^*: H^{\frac{1}{2}}(\partial\Omega) \to H^1(\Omega) \tag{8}$$

and

$$\chi_\Gamma: H^{\frac{1}{2}}(\partial\Omega) \to H^{\frac{1}{2}}(\Gamma) \tag{9}$$

$\chi_\Gamma^*$ is the opertor restriction to $\Gamma$ which is the adjoint $\chi_\Gamma$ given by

$$\chi_\Gamma^*: H^{\frac{1}{2}}(\Gamma) \to H^{\frac{1}{2}}(\partial\Omega) \tag{10}$$

**Definition 2.1.** The autonmous system asociated with system (1- 2) is said to be regionally boundary observable (or $\Gamma$- observable) if

$$\text{Im}(\chi_\Gamma \gamma_0 K^*) = H^{\frac{1}{2}}(\Gamma) \tag{11}$$

## 2.2 $\Gamma_E$- Observation

In this subsection an intrigued augmentation of territorial case as in ref. [11] is created to the limit case. Consequently, the portrayal of provincial limit exponential observation needs some notions which are related to the exponential behavior that are stability, detectability, and observer. we characterize the ideas which are identified with $\Gamma_E$-Observbility, what's more, give an essential hypothesis which gives an exponential onlooker for the first framework in basic subregion $\Gamma$.

The regnal bundary exponential observer in $\Gamma$ might be viewed as interior regional exponential observer in ω if we assume the following.

- Let $\Re$ be the continuous liner extension operator [15]

  $\Re : H^{\frac{1}{2}}(\partial\Omega) \longrightarrow H^1(\Omega)$ such that

  $$\gamma_0 \Re h(\mu, t) = h(\mu, t), \quad \forall h \in H^{\frac{1}{2}}(\partial\Omega) \tag{12}$$

- Let $r > 0$ is an arbitrary and sufficiently small real and let the sets

  $$E = \bigcup_{z \in \Gamma} B(z, r) \text{ and } \omega_r = E \cap \Omega \tag{13}$$

where $B(z, r)$ is the ball of range r focused in $z(\mu, t)$ and $\Gamma$ is a part of $\bar{\omega}_r$(fig.2).

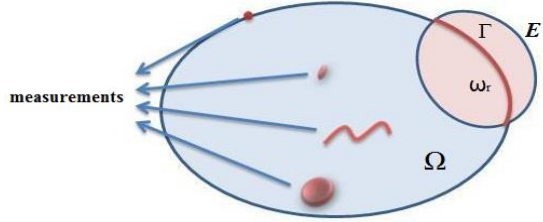

Fig. 2: The whole domain $\Omega$ and the set $E$.

**Definition 2.6.** The semi-group $(S_A(t))t \geq 0$ is regionally boundary exponentially stable *in* $H^{\frac{1}{2}}(\Gamma)$(or $\Gamma_E$-stable)if, for somepositive constants $F_\Gamma$ & $\sigma_\Gamma$, then

$$\|\chi_\Gamma \gamma_0 S_A(.)\|_{H^{\frac{1}{2}}(\Gamma)} \leq F_\Gamma e^{-\sigma_\Gamma t}, t \geq 0 \tag{14}$$

If $(S_A(t))_{t\geq 0}$ is $\Gamma_E$-stable, then for every $z_0(.) \in H^1(\overline{\Omega})$ the solution of autnmous system associated with (1) satisfies

$$\lim_{t\to\infty}\|z(t)\|_{H^{\frac{1}{2}}(\Gamma)} = \lim_{t\to\infty}\|\chi_\Gamma\gamma_0 S_A(t)z_0\|_{H^{\frac{1}{2}}(\Gamma)} = 0 \tag{15}$$

**Definition2.7.** The system (1) is regionally boundary exponentially stable on $\Gamma$ (or$\Gamma_E$-stable)if, the operator A genertes astrongly contnuous semi - group$(S_A(t))_{t\geq 0}$ which is $\Gamma_E$-stable.

**Remark2.8.** If the system (1) is $\Gamma_E$-stable, then solution of autonmous system with (1) convrges exponentially to zero when $t$ tends to $\infty$.

**Definition2.9.** The system (1) - (2) is locally limit exponentially noticeable on $\Gamma$ (or $\Gamma_E$-detectble) if there exst anoperator $H_\Gamma: \mathcal{O} \to H^{\frac{1}{2}}(\Gamma)$ such that the operator $A - H_\Gamma C$ produces a firmly constant semigroup $(S_{H_\Gamma}(t))_{t\geq 0}$ which is $\Gamma_E$-stable.

Now, assume that the dynamic framework:

$$\frac{\partial w}{\partial t}(\mu, t) = L_\Gamma w(\mu, t) + G_\Gamma u(t) + H_\Gamma y(t) \quad \Omega \times (0, \infty)$$

$$w(\mu, 0) = w_0(\mu) \quad\quad\quad\quad\quad\quad\quad\quad\quad\quad \overline{\Omega} \tag{16}$$

$$\frac{\partial w}{\partial v_{L_\Gamma}}(\mu, t) = 0 \quad\quad\quad\quad\quad\quad\quad\quad\quad\quad \partial\Omega \times (0, \infty)$$

where $L_\Gamma$ produces a firmly constant semigroup $(S_{L_\Gamma}(t))_{t\geq 0}$ which is stabile on the objec space W, $G_\Gamma \in L(R^p, W)$ and $H_\Gamma \in L(R^q, W)$. The system (16) defines an $\Gamma_E$-estimater for $\alpha_\Gamma z(\mu, t)$ if

(1) $\lim_{t\to\infty}\|w(.,t) - \alpha_\Gamma z(.,t)\|_{H^{\frac{1}{2}}(\Gamma)} = 0 \tag{17}$

(2) $\alpha_\Gamma$ maps D(A) into D($L_\Gamma$) where $\alpha_\Gamma = \chi_\Gamma\gamma_0 T$ and w(.,t) is the solution of system(16).

**Definition 2.10.** The system (16) is regionaly bundary exponential observer on $\Gamma$ (or $\Gamma_E$-observer) for the system (1)-(2) if the following conditions hold:

(1) The exist operators $M_\Gamma \in L\left(\mathcal{O}, H^{\frac{1}{2}}(\Gamma)\right)$ and $N_\Gamma \in L\left(H^{\frac{1}{2}}(\Gamma)\right)$ such that:

$$M_\Gamma C + N_\Gamma \alpha_\Gamma = I_\Gamma \tag{18}$$

$$\alpha_\Gamma A + L_\Gamma \alpha_\Gamma = H_\Gamma C \text{ and } G_\Gamma = \alpha_\Gamma B$$

(2) The system (16) defines an $\Gamma_E$-estimator for the state of original system.

**Definition 2.11.** The system (16) is regionaly bundary exponential identity observer (or identity $\Gamma_E$-observer) for the system (1)-(2) if Z = W and $\alpha_\Gamma = I_\Gamma$. In thiscase, wehave

$$L_\Gamma = A - H_\Gamma C \text{ and } G_\Gamma = B \tag{19}$$

And after that, the dynamical framework (16) moves toward becoming:

$$\frac{\partial w}{\partial t}(\mu, t) = Aw(\mu, t) + Bu(t) + H_\Gamma(Cw(\mu, t) - y(\mu, t)) \quad \Omega \times (0, \infty)$$

$$w(\mu, 0) = 0 \qquad\qquad\qquad\qquad \bar{\Omega} \tag{20}$$

$$\frac{\partial w}{\partial v_A}(\mu, t) = 0 \qquad\qquad\qquad\qquad \partial\Omega \times (0, \infty)$$

**Definition 2.12.** We can say that the system (1)-(2) is regionally boundary exponentially observable on $\Gamma$ (or $\Gamma_E$-observable) if, there exist a dynamical system (16) which is $\Gamma_E$-observer for the original system.

## 3. Regional boundary exponentially observation error

The general observation problem in deciding whether the knowledge of the output together with the system dynamics makes the state reconstruction possible. When the system is exponentially observable, the state reconstruction leads necessarily to an exponentially reconstruction error, also called the exponentially observation error. In the usual observation problem, we suppose that error $Er(z)$ defined by:

$$Er(z) = \|z_0(.) - \tilde{z}_0(.)\|^2_{L^2(0,T;\mathbb{R}^q)}, \tag{21}$$

where $\tilde{z}_0(.)$ holds for the exponentially reconstrction state and $z_0(.)$ for the original state.

In our case, the exponentially observation error unmistakably relies upon the objective locale $\Gamma$ where the state is to be observed together with the structure and number of sensors and the exponentially reconstruction state $\tilde{z}_0(.)$ is the state of $\Gamma_E$- observer system (16) which is $w_0(.)$:

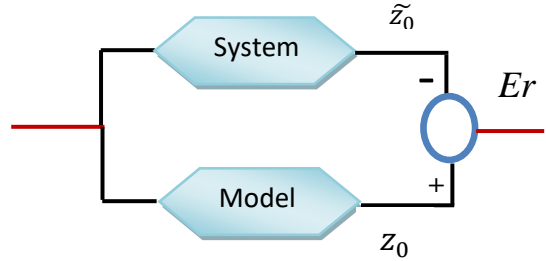

**Fig. 3. the (exponential) ly observation error**

Now, let we consider the following sets:

The set of exponentially observable states in $\Gamma$ is denoted by

$$\Psi_{\Gamma_E} = \{\Gamma_E\text{-observable states}\} \tag{22}$$

The set $\Psi_{\Omega_E}$ stands for the set of exponentially observable states (in $\Omega$). The associated exponentially observation error is then denoted by

$$Er_{\Gamma_E} : z_0 \in \Psi_{\Gamma_E} \rightarrow Er_{\Gamma_E}(z_0) \in \mathbb{R} \tag{23}$$

The set of exponentially observable states in $\Gamma$ at the point when estimations are acquired by methods for one zone sensor $C$ is denoted by

$$\Psi^C_{\Gamma_E} = \{\Gamma_E\text{-observable states by means of the sensor } C\} \tag{24}$$

The exponentially observation error associated with the set (24) is denoted by

$$Er_{\Gamma_E}(\cdot, C): z_0 \in \Psi_{\Gamma_E}^C \longrightarrow Er_{\Gamma_E}(z_0, C) \in \mathbb{R} \tag{25}$$

When measurements are obtained by $k$ sensors $(C_i) 1 \leq i \leq k$, we consider the set

$$\Psi_{\Gamma_E}^{C_{\{1,2,\dots,k\}}} = \{\Gamma_E\text{-observable states by means of the sensor } C_1, C_2, \dots, C_k\} \tag{26}$$

The associated exponentially observation error is then denoted by

$$Er_{\Gamma_E}(\cdot,): z_0 \in \Psi_{\Gamma_E}^{C_{\{1,2,\dots,k\}}} \longrightarrow Er_{\Gamma_E}(z_0, C_1, C_2, \dots, C_k) \in \mathbb{R} \tag{27}$$

**Remark 3.1:**

For any nonempty boundary region $\Gamma$, $\Gamma \subset \partial\Omega$, we have

$$\Psi_{\Omega_E} \subset \Psi_{\Gamma_E} \tag{28}$$

That's mean the exponential observability implies regional boundary exponential observability.

**Remark 3.2:**

(i) on account of pointwise sensors, in the above notation the sensor $C$ is replaced by the location $b$. The sets given by (24) and (26) become respectively $\Psi_{\Gamma_E}^b$ and $\Psi_{\Gamma_E}^{b_{\{1,2,\dots,k\}}}$.

(ii) On account of limit zone (individually pointwise) sensors, the documentation continues as before aside from that the sensor underpins (resp. areas) are subsats of the lmit $\partial\Omega$ of the domain $\Omega$.

The exponantially reconstrction method is based on finding a state $z_0^*$ which realizes

$$Er(z_0^*) = \min_{z_0} Er(z_0) \tag{29}$$

**Note that** the minimum $Er(z_0^*)$ relies upon how the estimations have been considred, i.e., the number and structure of controls. In this manner we compose

$$Er(z_0^*) = Er(z_0^*, C) \tag{30}$$

for the exponentially observation error depending on the sensor $C$. Furthrmore, we supose that the mapping $C \longrightarrow Er(z_0^*, C)$ is continuous. This means that if the sensor parameters are slightly modified, then the exponentially observation error changes continuously. on account of $k$ sensors, i.e., a multi-yield circumstance, it will be meant by

$$Er(z_0^*) = Er(z_0^*, C_1, C_2, \dots, C_k) \tag{31}$$

When the controls are pointwise, we write

$$Er(z_0^*) = Er(z_0^*, b_1, b_2, \ldots, b_k) \qquad (32)$$

The notation remains similar in the boundary case except that the supports and locations of sensors are subsets of $\partial\Omega$.

In what follows we consider the observtion error function defined by (23), (25) and (27) depending on the case considred. We asume that the functions $Er, Er(\cdot, C)$ and $Er(\cdot, C_1, C_2, \ldots, C_k)$ are convex with respect to the initial state $z_0$. Therefore we immediately have the follwing important result.

**Proposition 3.3:**

*Let $\Gamma$ be a given region of $\partial\Omega$. Suppose that the system (1)-(2) is $\Gamma_E$- observable. Then the exponentially observation error in $\Gamma$ is lower than that in the whole domain $\Omega$.*

*Proof.*

Consider the sets $\Psi_{\Gamma_E}$ and $\Psi_{\Omega_E}$. Since $\Gamma \subset \partial\Omega \subset \Omega$, by Remark 3.2, any exponentially observable state in $\Omega$ is $\Gamma_E$- *observable* and then  $\Psi_{\Omega_E} \subset \Psi_{\Gamma_E}$ , because exponential observability implies regional boundary exponential observability for any $\Gamma \subset \Omega$. Therefore,

$$\min_{z_0 \in \Psi_{\Gamma_E}} Er(z_0) \leq \min_{z_0 \in \Psi_{\Omega_E}} Er(z_0) \qquad \blacksquare \qquad (33)$$

**Remark 3.4:**

(i) As a speculation of the above outcome, we can exhibit that if $\Gamma_1$ and $\Gamma_2$ are two subdomans of $\partial\Omega$ such that $\Gamma_1 \subset \Gamma_2$, then

$$\mathop{\text{Min}}_{z_0 \in \Psi_{\Gamma_{1E}}} Er(z_0) \leq \min_{z_0 \in \Psi_{\Gamma_{2E}}} Er(z_0) \qquad (34)$$

The prof is similar to above result and proceeds by replacing $\Gamma_E$ by $\Gamma_{1E}$ and $\Omega_E$ by $\Gamma_{2E}$.

(ii) The above result means that the smaller the boundary region $\Gamma$ where the state is to be exponential observed, the lower the corrsponding exponential observtion error. The result also true even if the region considered is internal region on $\Omega$, with usual observability[2].

The exponential observation error on $\Gamma$ is assumed to be convex, and we write

$$Er(z_0^*) = \min_{z_0 \in \Psi_{\Gamma_E}} Er(z_0) \tag{35}$$

And

$$Er_\Omega(z_0^*) = Er(z_0^*) \; for \min_{z_0 \in \Psi_{\Omega_E}} Er(z_0)$$

(36)

## 4. Conclusion

  The ideas created in this paper is identified with the exponentially (perception) error in limit locale. By presenting the thought of exponential recognizability, we demonstrated a natural outcome which stipulates that the exponentially perception mistake in limit locale is lower than in the entire space. Many inquiries still opened. This is the situation of, for instance, the issue of finding the ideal sensor area guaranteeing such a goal.